\documentclass[3p,10pt,a4paper,twoside,fleqn,sort&compress]{filomat}
\usepackage{amssymb,amsmath,latexsym}
\usepackage[varg]{pxfonts}

\newtheorem{theorem}{Theorem}
\newtheorem{corollary}{Corollary}
\newtheorem{remark}{Remark}

\def\const{\mathop{\rm const}}
\def\Ai{\mathop{\rm Ai}\nolimits}
\def\Bi{\mathop{\rm Bi}\nolimits}

\newcommand{\FilVol}{xx}
\newcommand{\FilYear}{yyyy}
\newcommand{\FilPages}{zzz--zzz}
\newcommand{\FilDOI}{(will be added later)}
\newcommand{\FilReceived}{Received: dd Month yyyy; Accepted: dd Month yyyy}

\begin{document}

\title{A note about the torsion of null curves in the $3$-dimensional Minkowski spacetime and the Schwarzian derivative}

\author{Zbigniew Olszak}
\ead{zbigniew.olszak@pwr.edu.pl}
\address{Wroc{\l}aw University of Technology, Institute of Mathematics and Computer Science, 50-370 Wroc{\l}aw, Poland}
\newcommand{\AuthorNames}{Z. Olszak}

\newcommand{\FilMSC}{Primary 53A35; Secondary 53B30, 53B50, 53C50.}
\newcommand{\FilKeywords}{Minkowski spacetime of dimension 3; null curve; light-like curve; torsion of a curve; null helix; slant helix; Airy function; Schwarzian derivative.}
\newcommand{\FilCommunicated}{(name of the Editor, mandatory)}

\begin{abstract}
The main topic of this paper is to show that in the 3-dimensional Minkowski spacetime, the torsion of a null curve is equal to the Schwarzian derivative of a certain function appearing in a description of the curve. As applications, we obtain descriptions of the slant helices, and null curves for which the torsion is of the form $\tau=-2\lambda s$, $s$ being the pseudo-arc parameter and $\lambda=\const\neq0$.
\end{abstract}

\maketitle

\makeatletter
\renewcommand\@makefnmark%
{\mbox{\textsuperscript{\normalfont\@thefnmark)}}}
\makeatother

\section{Introduction}

There are very many papers about geometric properties of null curves in the Minkowski spacetimes. We refer the monographs \cite{DB, DJ}, and the survey articles \cite{Dug, IL, L}, etc.

On the other hand, there is the classical notion of the Schwarzian derivative in mathematical analysis. This notion has many important applications in mathematical analysis (real and complex) and differential geometry; see \cite{DG, DO, O, OS, OT}, etc. The author is specially inspired by the paper \cite{DO}, where it is shown a strict relation between the Schwarzian derivative and the curvature of worldlines in 2-dimenional Lorentzian manifolds of constant curvature.

In the presented short paper, we will show that the torsion of a null curve in the 3-dimensional Minkowski spacetime $\mathbb E^3_1$ is equal to the Schwarzian derivative of a certain function appearing in a description of the curve. Descriptions of the slant helices are obtained, and null curves for which the torsion is given by $\tau=-2\lambda s$, $s$ being the pseudo-arc parameter and $\lambda=\const\neq0$.

\section{Preliminaries}

Let  $\mathbb E^3_1$ be the 3-dimensional Minkowski spacetime, that is, the Cartesian $\mathbb R^3$ endowed with the standard Minkowski metric $g$ given with respect to the Cartesian coordinates $(x,y,z)$ by
\begin{equation}
\label{mink1}
  g = dx\otimes dx + dy\otimes dy - dz\otimes dz,
\end{equation}
or as the symmetric 2-form $g = dx^2 + dy^2 - dz^2$. 

Let $\alpha\colon I\to\mathbb E^3_1$ be a null (ligth-like) curve in $\mathbb E^3_1$, $I$ being an open interval. Thus, $g(\alpha^{\,\prime},\alpha^{\,\prime})=0$, that is, $g(\alpha^{\,\prime}(t),\alpha^{\,\prime}(t))=0$ for any $t\in I$. We also assume that the curve is non-degenerate, in the sense the three vector fields $\alpha^{\,\prime}$, $\alpha^{\,\prime\prime}$, $\alpha^{\,\prime\prime\prime}$ are linearly independent at every point of the curve. 

Since $g(\alpha^{\,\prime},\alpha^{\,\prime})=0$ and $g(\alpha^{\,\prime},\alpha^{\,\prime\prime})=0$, it must be that $g(\alpha^{\,\prime\prime},\alpha^{\,\prime\prime})>0$. A parametrization of the null curve is said to be pseudo-arc (or distinguished) if $g(\alpha^{\,\prime\prime},\alpha^{\,\prime\prime})=1$. A null curve can always be parametrized by a pseudo-arc parameter. However, such a parameter is not uniquely defined. Precisely, for a null curve $\alpha$, if $s_1$ is a pseudo-arc parameter, then $s_2$ is a pseudo-arc parameter if and only if there exists a constant $c$ such that $s_2=\pm s_1+c$. 

In the sequel, we assume that the parametrization of a null curve is pseudo-arc, and we denote such a parameter by $s$.

In the next section, we need the standard theorms concerning of null curves which can be formulated in the following manner (see e.g. \cite{DJ, Dug, IL, L}):

Let  $\alpha$ be a null curve in the 3-dimensional Minkowski spacetime $\mathbb E^3_1$. Then, there exists the only one Cartan moving frame $(\mathbf L=\alpha^{\,\prime}, \mathbf N, \mathbf W)$ and the function $\tau$ defined along the curve $\alpha$ and such that 
\begin{equation}
\label{scpr}
  g(\mathbf L,\mathbf N) = 
  g(\mathbf W,\mathbf W) = 1,\quad
  g(\mathbf L,\mathbf L) = 
  g(\mathbf L,\mathbf W) = 
  g(\mathbf N,\mathbf N) = 
  g(\mathbf N,\mathbf W) = 0, 
\end{equation}
and the following system of differential equations 
\begin{equation}
\label{CFequ}
  \mathbf L^{\,\prime} = \mathbf W, \quad
  \mathbf N^{\,\prime} = \tau \mathbf W, \quad
  \mathbf W^{\,\prime} =\null- \tau \mathbf L - \mathbf N
\end{equation}
is satisfied. These vector fileds are given by
\begin{equation}
\label{M1}
  \mathbf L = \alpha^{\,\prime}, \quad
  \mathbf W = \alpha^{\,\prime\prime}, \quad
  \mathbf N = \null- \alpha^{\,\prime\prime\prime} 
								- \frac12 g(\alpha^{\,\prime\prime\prime},\alpha^{\,\prime\prime\prime})
								\alpha^{\,\prime},
\end{equation}
and the function $\tau$ by 
\begin{equation}
\label{tau}
  \tau = \frac12 g(\alpha^{\,\prime\prime\prime},\alpha^{\,\prime\prime\prime}).
\end{equation}

From these results it can be deduced that a given function $\tau$ on an open interval $I$, there exists the only one null curve $\alpha\colon I\to\mathbb E^3_1$ realizing (\ref{scpr}) and (\ref{CFequ}) up to the orientation of this curve and up to the isometries of the Minkowski space $\mathbb E^3_1$. 

The triple $(\mathbf L, \mathbf N, \mathbf W)$ defined in (\ref{M1}) is called the Frenet frame, the function $\tau$ defined in (\ref{tau}) is called the torsion, and the equations (\ref{CFequ}) are called the Frenet equations of the null curve $\alpha$. Since
\begin{equation}
\label{orient}
  \det[\mathbf L, \mathbf N, \mathbf W] 
		= \det[\alpha^{\,\prime}, \alpha^{\,\prime\prime}, \alpha^{\,\prime\prime\prime}],
\end{equation}
the frames $(\mathbf L, \mathbf N, \mathbf W)$ and $\left(\alpha^{\,\prime}, \alpha^{\,\prime\prime}, \alpha^{\,\prime\prime\prime}\right)$ have the same orientations.

In the following section, we are going to expresse the torsion $\tau$ and the frame $(\mathbf L,\mathbf N,\mathbf W)$ with the help of a special function related to a pseudo-arc parametrization of a null curve in $\mathbf E^3_1$.

\section{A description of the torsion}

Let $\alpha\colon I\to\mathbb E^3_1$ be a null curve. Simplifying denotations, we write 
$\alpha(s) = (x(s),y(s),z(s))$, $s\in I$, where $s$ is a pseudo-arc parameter, and $x(s)$, $y(s)$, $z(s)$ are certain functions of $s$. Then, we have 
\begin{equation*}
  \alpha^{\,\prime} = x^{\,\prime} \dfrac{\partial}{\partial x}\Big|_{\alpha}
		 + y^{\,\prime} \dfrac{\partial}{\partial y}\Big|_{\alpha}
			+ z^{\,\prime} \dfrac{\partial}{\partial z}\Big|_{\alpha}.
\end{equation*}
For simplicity, instead of that, we will write $\alpha^{\,\prime} = (x^{\,\prime},y^{\,\prime},z^{\,\prime})$. And in the similar manner, the next derivatives of $\alpha$ will be written, e.g., $\alpha^{\,\prime\prime} = (x^{\,\prime\prime},y^{\,\prime\prime},z^{\,\prime\prime})$.

Using (\ref{mink1}), our two assumptions: $g(\alpha^{\,\prime},\alpha^{\,\prime})=0$ (the nullity condition), and $g(\alpha^{\,\prime\prime},\alpha^{\,\prime\prime})=1$ (the pseudo-arc parametrization) give the following two equalities
\begin{eqnarray}
\label{ds:equ1}
  x^{\,\prime\,2} + y^{\,\prime\,2} - z^{\,\prime\,2} &=& 0, \\
\label{ds:equ2}
  x^{\,\prime\prime\,2} + y^{\,\prime\prime\,2} - z^{\,\prime\prime\,2} &=& 1.
\end{eqnarray}
One notes that the shapes of the equalities (\ref{ds:equ1}) and (\ref{ds:equ2}) exclude the situation when at least one of the functions $x^{\,\prime}$, $y^{\,\prime}$, $z^{\,\prime}$ vanishes on an open subinterval of $I$. In the sequel, restricting slightly the assumptions, we will consider only the case when $x^{\,\prime}\neq0$, $y^{\,\prime}\neq0$ and $z^{\,\prime}\neq0$ on $I$. 

It is a standard and elementary idea that from (\ref{ds:equ1}), it follows that 
\begin{equation}
\label{fh}
		x^{\,\prime} = h,\quad 
  y^{\,\prime} = \frac{h}{2} \left(f - \frac{1}{f}\right),\quad
  z^{\,\prime} = \frac{h}{2} \left(f + \frac{1}{f}\right),
\end{equation}
$f$ and $h$ being certain non-zero functions on $I$. Hence,
\begin{eqnarray*}
		x^{\,\prime\prime} &=& h^{\prime}, \\
  y^{\,\prime\prime} &=& \frac{fh^{\prime}(f^2 - 1) + h f^{\prime}(f^2 + 1)}{2f^2}, \\
  z^{\,\prime\prime} &=& \frac{fh^{\prime}(f^2 + 1) + h f^{\prime}(f^2 - 1)}{2f^2}. 
\end{eqnarray*}
In view of the above relations, the equality (\ref{ds:equ2}) turns into $h^2f^{\,\prime\,2}=f^2$. Hence, $f^{\,\prime}$ is non-zero (and has constant sign) on $I$. Consequently, 
\begin{equation*}
  h = \varepsilon \frac{f}{f^{\,\prime}}, \ \varepsilon=\pm1.
\end{equation*}

Thus, for the vector field $\mathbf L$ (cf.\ (\ref{M1})), we have 
\begin{equation}
\label{ds:L}
  \mathbf L 
		= \alpha^{\,\prime} 
		= \frac{\varepsilon}{2f^{\,\prime}} \left(2f,f^2 - 1,f^2 + 1\right).
\end{equation}

Consequently, we get the following description of the curve $\alpha$ 
\begin{equation*}
  \alpha(s) = \alpha(s_0)
     + \frac{\varepsilon}{2} \int_{s_0}^s \frac{1}{f^{\,\prime}(t)}
		    \left(2f(t),f^2(t) - 1,f^2(t) + 1\right)\,dt,\ s,s_0\in I.
\end{equation*}

Conversely, if a curve $\alpha$ is given by the last formula, then (\ref{ds:equ1}) and (\ref{ds:equ2}) are fulfilled so that the curve is null and not geodesic, and the parameter $s$ is distinguish. 

From (\ref{ds:L}), we obtain for the vector field $\mathbf W$ (cf.\ (\ref{M1})), 
\begin{equation}
\label{ds:W}
  \mathbf W = \alpha^{\,\prime\prime} 
    = \null-\frac{\varepsilon f^{\,\prime\prime}}{2f^{\,\prime\,2}}
				    \left(2f,f^2-1,f^2+1\right) + \varepsilon \left(1,f,f\right).
\end{equation}

From (\ref{ds:W}), we find 
\begin{equation}
\label{ds:tres}
  \alpha^{\,\prime\prime\prime}
    = \varepsilon \frac{2f^{\,\prime\prime\,2}-f^{\,\prime}f^{\,\prime\prime\prime}}
        {2f^{\,\prime\,3}} \left(2f,f^2 - 1,f^2 + 1\right) 
        - \frac{\varepsilon f^{\,\prime\prime}}{f^{\,\prime}} (1,f,f)
        + \varepsilon f^{\,\prime}(0,1,1). 
\end{equation}

To compute $g(\alpha^{\,\prime\prime\prime},\alpha^{\,\prime\prime\prime})$, using (\ref{mink1}), we find at first the following 
\begin{eqnarray*}
  && g\left(\left(2f,f^2 - 1,f^2 + 1\right),\left(2f,f^2 - 1,f^2 + 1\right)\right) = 0, \\
		&& g\left(\left(2f,f^2 - 1,f^2 + 1\right),(1,f,f)\right) = 0, \quad
		   g\left(\left(2f,f^2 - 1,f^2 + 1\right),(0,1,1)\right) = -2, \\
		&& g\left((1,f,f),(1,f,f)\right) = 1, \quad 
		   g\left((1,f,f),(0,1,1)\right) = 0, \quad
		   g\left((0,1,1),(0,1,1)\right) = 0. 
\end{eqnarray*}
Then, having (\ref{ds:tres}) and applying the above formulas, we get 
\begin{equation}
\label{ds:gtres}
  g(\alpha^{\,\prime\prime\prime},\alpha^{\,\prime\prime\prime}) 
		= \frac{2f^{\,\prime}f^{\,\prime\prime\prime}- 3f^{\,\prime\prime\,2}}{f^{\,\prime\,2}}. 
\end{equation}

In view of (\ref{ds:gtres}) and (\ref{tau}), the torsion must be of the form 
\begin{equation}
\label{ds:tau}
  \tau =\frac{2f^{\,\prime}f^{\,\prime\prime\prime}- 3f^{\,\prime\prime\,2}}
		              {2f^{\,\prime\,2}}
							= \left(\frac{f^{\,\prime\prime}}{f^{\,\prime}}\right)^{\,\prime}
           - \frac12 \left(\frac{f^{\,\prime\prime}}{f^{\,\prime}}\right)^2.
\end{equation}
Now, it is important to note that the right hand side of the formula (\ref{ds:tau}) is just the Schwarzian derivative of the function $f$, which is usually denoted by $S(f)$. Thus, $\tau=S(f)$.  

Finally, applying (\ref{ds:L}), (\ref{ds:tres}) and (\ref{ds:gtres}) into (\ref{M1}), we find the vector field 
\begin{eqnarray}
\label{ds:N}
  \quad \mathbf N 
		= \null- \frac{\varepsilon f^{\,\prime\prime\,2}}{4f^{\,\prime\,3}}	
		   \left(2f,f^2 - 1,f^2 + 1\right) 
     + \frac{\varepsilon f^{\,\prime\prime}}{f^{\,\prime}} (1,f,f)
     - \varepsilon f^{\,\prime} (0,1,1). 
\end{eqnarray}

Summarizing the above considerations, we can formulate the following theorem. 

\begin{theorem}
Let $\mathbb E^3_1$ be the 3-dimensional Minkowski spacetime. Any (non-degenerate) null curve $\alpha$ in $\mathbb E^3_1$ can be parametrized in the following way 
\begin{equation}
\label{ds:curv}
  \alpha(s) = \alpha(s_0)
     + \frac{\varepsilon}{2} \int_{s_0}^s \frac{1}{f^{\,\prime}(t)}
		    \left(2f(t),f^2(t) - 1,f^2(t) + 1\right)\,dt,\ s,s_0\in I,
\end{equation}
where $s$ is a pseudo-arc parameter, $I$ is a certain open interval, $f$ is a non-zero function with non-zero derivative $f^{\,\prime}$ on $I$. The torsion $\tau$ of such a curve is equal to the Schwarzian derivative of the function $f$, that is,
\begin{equation}
\label{ds:tau2}
  \tau = S(f) = \left(\frac{f^{\,\prime\prime}}{f^{\,\prime}}\right)^{\,\prime}
           - \frac12 \left(\frac{f^{\,\prime\prime}}{f^{\,\prime}}\right)^2. 
\end{equation}
The vector fields forming the Frenet frame of the curve $\alpha$ are given by the formulas (\ref{ds:L}), (\ref{ds:W}) and (\ref{ds:N}). 
\end{theorem}

\begin{remark} 
Applying formulas (\ref{ds:L}), (\ref{ds:W}), (\ref{ds:N}), it can be verified that 
\begin{equation*}
  \det[\mathbf L, \mathbf N, \mathbf W] = \varepsilon.
\end{equation*}
This together with (\ref{orient}) implies that the constant $\varepsilon$ appering in (\ref{ds:curv}) corresponds to the orientation of the curve $\alpha$. Note that the torsion does not depend on the orientation of the curve. Moreover, the torsion and the orientation does not depend on the sign of the function $f$.
\end{remark}

\begin{remark}
The Schwarzian derivative $S$ is an invariant of a fractional-linear transformation $T$ of the 1-dimensional real projective space $\mathbb RP^1=\mathbb R\cup\infty$ (cf.\ e.g.\ \cite{O}). That is, 
\begin{math}
  S(T\circ f)=S(f)
\end{math}
if $f$ is a function on $\mathbb RP^1$ and $T\colon\mathbb RP^1\to\mathbb RP^1$ is given by 
\begin{equation}
  T(r) = \frac{ar+b}{cr+d},\ r\in\mathbb RP^1,\ a,b,c,d\in\mathbb R, ad-bc\neq0.
\end{equation}

We can apply the above fact seeking for null curves with given torsion $\tau$. However, we should be careful since the domains of our functions $f$ and $T\circ f$ may be defined only on some open subintervals lying on the real line $\mathbb R$. 
\end{remark}

\section{Null Cartan helices}

It is well-known that there are exactly three types of null curves with constant torsion in the Minkowski spacetime $\mathbb E^3_1$ (cf e.g., \cite{FGL2}) up to the orientation of the curve and up to the isometries of the space. They are often called the null Cartan helices. 

As a first application of the results from the previous section, we demonstrate how these classes of curves can be recovered from their torsions.

(a) For $f(s)=s$, it holds $S(f)=0$. In (\ref{ds:curv}), we put
\begin{math}
  f(s)=s, \ s_0=0, \ \alpha(s_0)=(0,0,0), \ \varepsilon=1.
\end{math}
Then, we obtain the curve 
\begin{equation*}
  \alpha(s) = \frac{1}{6}\left(3s^2, s^3 - 3s, s^3 + 3s\right), 
\end{equation*}
for which by (\ref{ds:tau2}) we have $\tau=0$. Thus, $\alpha$ is a positively oriented null Cartan helix of zero torsion. 

(b) For $f(s)=-\cot(cs/2)$, it holds $S(f)=c^2/2$. In (\ref{ds:curv}), we put
\begin{equation*}
  f(s)=-\cot\frac{cs}{2},  
		\ \alpha(0)=\left(\frac{1}{c^2},0,0\right), 
		\ \varepsilon=1, 
		\ c=\const.>0.
\end{equation*}
Then, we obtain the curve 
\begin{equation*}
  \alpha(s) = \frac{1}{c^2}\left(\cos(cs),\sin(cs),cs\right),
\end{equation*}
for which by (\ref{ds:tau2}) it holds $\tau=c^2/2$. Thus, $\alpha$ is a positively oriented null Cartan helix of constant positive torsion. 

(c) For $f(s)=e^{cs}$, it holds $S(f)=-c^2/2$. In (\ref{ds:curv}), we put
\begin{equation*}
  f(s)=e^{cs}, 
		\ \alpha(0)= \left(0,\frac{1}{c^2},0\right), 
		\ \varepsilon=1,
		\ c=\const.>0.
\end{equation*}
Then, we obtain the curve 
\begin{equation*}
  \alpha(s) = \frac{1}{c^2}\left(cs, \cosh(cs), \sinh(cs)\right),
\end{equation*}
for which by (\ref{ds:tau2}) we have $\tau=-c^2/2$. Thus, $\alpha$ is a positively oriented null Cartan helix of constant negative torsion. 

Thus, we have seen the following:

\begin{corollary}
Null helices in $\mathbb E^3_1$ form the three classes described in (a) -- (c) in the above. The description is valid up to the pseudo-arc parameter changies, up to the orientation of the curve, and up to the isometries of the space. 
\end{corollary}

A curve $\alpha\colon I\to\mathbb E^3_1$ is called a general (or generalized) helix if there exists a non-zero vector $V$ in $\mathbb E^3_1$ such that $g(\alpha^{\,\prime},V)=\const.$; cf. \cite{FGL1,FGL2,SKG}, etc. This means that tangent indicatrix is laid in a plane or, equivalently, there exists a non-zero constant vector $V$ in $\mathbb E^3_1$ for which $g(\alpha^{\,\prime\prime},V)=0$, that is, $V$ is orthogonal to the acceleration vector field $\alpha^{\,\prime\prime}$. 

For null curves, it is already proved that null general helices in $\mathbb E^3_1$ are precisely the null Cartan helices; cf. ibidem. 

\section{Null slant helices}

Following the ideas of \cite{AL,CK,GG}, a slant helix is defined to be the curve (null as well as non-null) in $\mathbb E^3_1$ which satisfies the condition 
\begin{equation}
\label{slh}
  g(\alpha^{\,\prime\prime},V) = c = \const.
\end{equation}
along the curve $\alpha$, where $V$ is a constant vector. Thus, a general helix is a slant helix with $c=0$. Conversely, a slant helix with $c=0$ becomes a general helix. In \cite[Theorem 1.4]{AL}, it is proved that a null curve in $\mathbb E^3_1$ is a slant helix if and only if its torsion is given by 
\begin{equation}
\label{taush}
  \tau = \frac{a}{(cs+b)^2}, \  a,b,c= \const.,
\end{equation}
where $c$ is just the constant realizing (\ref{slh}). 

As the second applications of the results from Section 3, we will describe the null slant helices in $\mathbb E^3_1$ which are different from the usual helices ($a\neq0$ and $c\neq0$ in (\ref{taush})). 

Note that moving the pseudo-arc parameter $s$ into $s-b/c$ and next modifying slightly the constant $a$, we can write the condition (\ref{taush}) as
\begin{equation}
\label{taush-2}
  \tau = \frac{a}{2s^2}, \  a= \const\neq0.
\end{equation}
We can also assume that $s>0$. Using (\ref{scpr}) and (\ref{CFequ}), it can be checked that when the relation (\ref{taush-2}) is fulfilled, then for the vector
\begin{equation*}
  V = \null- \frac{a}{2s} \mathbf L + s \mathbf N + \mathbf W
\end{equation*}
it holds $V^{\,\prime}=0$ and $g(\alpha^{\,\prime\prime},V) = g(\mathbf W,V) = 1$ (cf.\ ibidem). 

(a) In (\ref{ds:curv}), we put
\begin{equation*}
  f(s)=\ln s, 
		\ s_0=1, 
		\ \alpha(s_0)=\frac{1}{8}(-2,-1,3), 
		\ \varepsilon=1.
\end{equation*}
Then, we obtain the curve 
\begin{equation*}
  \alpha(s) = \frac{s^2}{8}\left(2(2\ln s-1),2\ln^2 s-2\ln s-1,2\ln^2 s-2\ln s+3\right), 
\end{equation*}
for which by (\ref{ds:tau2}) it holds 
\begin{equation*}
  \tau=S(f)=\frac{1}{2s^2}. 
\end{equation*}
Thus, $\alpha$ is a slant helix realizing (\ref{taush-2}) with $a=1$. 

(b) Let $a>1$ and $b=\sqrt{a-1}>0$. In (\ref{ds:curv}), we put
\begin{equation*}
  f(s)=\tan\left(\frac{1}{2}\ln s^b\right), 
		\ s_0=1, 
		\ \alpha(s_0)=\frac{1}{b}\left(-\frac{b}{b^2+4},-\frac{2}{b^2+4},\frac{1}{2}\right), 
		\ \varepsilon=1.
\end{equation*}
Then, we obtain the curve 
\begin{equation*}
  \alpha(s) = \frac{s^2}{b} 
		\left(\frac{2\sin(\ln s^b)-b\cos(\ln s^b)}{b^2+4},
							         -\frac{2\cos(\ln s^b)+b\sin(\ln s^b)}{b^2+4},
									    			\frac{1}{2}\right), 
\end{equation*}
for which by (\ref{ds:tau2}) it holds 
\begin{equation*}
  \tau=S(f)=\frac{1+b^2}{2s^2}=\frac{a}{2s^2}. 
\end{equation*}
Thus, $\alpha$ is a slant helix realizing (\ref{taush-2}) with $a>1$. 

(c) Let $0\neq a<1$. Then for $b=\sqrt{1-a}$, we have $b>0$ and $b\neq1$. Consider the case $a\neq-3$, that is, $b\neq2$. In (\ref{ds:curv}), we put
\begin{equation*}
  f(s)=s^{-b}, 
		\ s_0=1, 
		\ \alpha(s_0)=\frac{1}{2b}\left(-1,\frac{2b}{b^2-4},\frac{4}{b^2-4}\right), 
		\ \varepsilon=1.
\end{equation*}
Then, we obtain the curve 
\begin{equation*}
  \alpha(s) = \frac{s^2}{2b} 
		\left(-1,\frac{s^{-b}}{b-2}+\frac{s^b}{b+2},\frac{s^{-b}}{b-2}-\frac{s^b}{b+2}\right), 
\end{equation*}
for which by (\ref{ds:tau2}) it holds 
\begin{equation*}
  \tau=S(f)=\frac{1-b^2}{2s^2}=\frac{a}{2s^2}. 
\end{equation*}
Thus, $\alpha$ is a slant helix realizing (\ref{taush-2}) with $-3\neq a<1$. 

(d) In (\ref{ds:curv}), we put
\begin{equation*}
  f(s)=\frac{1}{s^2}, 
		\ s_0=1, 
		\ \alpha(s_0)=\frac{1}{16}(-4,1,-1), 
		\ \varepsilon=1.
\end{equation*}
Then, we obtain the curve 
\begin{equation*}
  \alpha(s) = \frac{1}{16}\left(-4s^2,s^4-4\ln s,\null-s^4-4\ln s\right), 
\end{equation*}
for which by (\ref{ds:tau2}) it holds 
\begin{equation*}
  \tau=S(f)=-\frac{3}{2s^2}. 
\end{equation*}
Thus, $\alpha$ is a slant helix realizing (\ref{taush-2}) with $a=-3$. 

Thus, we have shown the following:

\begin{corollary}
Null slant helices in $\mathbb E^3_1$ form the four classes described in (a) -- (d) in the above. The description is valid up to the pseudo-arc parameter changies, up to the orientation of the curve, and up to the isometries of the space. 
\end{corollary}

\section{Null curves with the torsion proportional to the pseudo-arc parameter}

In this section, we determine the null curves in $\mathbb E^3_1$ for which $\tau=-2\lambda s$, $\lambda=\const.\neq0$. We will use the formula (\ref{ds:tau2}). 

According to our Theorem, we need at first to find a solution of the differential equation 
\begin{equation}
\label{nEs}
  \left(\frac{f^{\,\prime\prime}}{f^{\,\prime}}\right)^{\,\prime}
           - \frac12 \left(\frac{f^{\,\prime\prime}}{f^{\,\prime}}\right)^2 = -2\lambda s.
\end{equation}
We seek for solutions of this equation in the form
\begin{equation}
\label{phi2f}
  f(s) = \int\frac{ds}{\phi^2(s)},
\end{equation}
$\phi$ being an unknown fucntion. Then the equation (\ref{nEs}) becomes the following differential equation
\begin{equation}
\label{eq4phi}
  \phi^{\,\prime\prime} - \lambda s \phi = 0.
\end{equation}
The general solution of the above equation is  
\begin{equation*}
  \phi(s) = c_1 \Ai\left(\mu s\right) + c_2 \Bi\left(\mu s\right), 
		\ \mu = \sqrt[3]{\lambda},\ c_1, c_2 = \const.
\end{equation*}
where  $\Ai$ and $\Bi$ are the Airy functions of the first and second kind, respectively. For the solutions of (\ref{eq4phi}) and for the special Airy functions, we refer \cite{PZ}, \cite{VS}, etc. In the below calculations, we use the basic properties of these functions. 

For our purpose, we take the only one solution of (\ref{eq4phi}), say 
\begin{math}
  \phi(s) = \Ai\left(\mu s\right). 
\end{math}
Then, from (\ref{phi2f}) we get 
\begin{equation}
\label{fnEs}
  f(s) = \frac{\pi}{\mu} \cdot 
		       \frac{\Bi(\mu s)}{\Ai(\mu s)}. 
\end{equation}
Next, 
\begin{equation}
\label{dnEs}
  f^{\,\prime}(s) = \frac{1}{\Ai^{\,2}(\mu s)}. 
\end{equation}
Having (\ref{ds:L}) with $\varepsilon=1$, and using (\ref{fnEs}) and (\ref{dnEs}), we can write $\alpha^{\,\prime}$ as
\begin{equation*}
  \alpha^{\,\prime}(s) 
		= \left( \frac{\pi}{\mu} \Ai(\mu s) \Bi(\mu s), 
	          \frac{1}{2\mu^2}
										 \left(\pi^2 \Bi^2(\mu s) - \mu^2 \Ai^2(\mu s)\right), \\
		         \frac{1}{2\mu^2}
									 	\left(\pi^2 \Bi^2(\mu s) + \mu^2 \Ai^2(\mu s)\right)\right). 
\end{equation*}
The integration of the last equality gives the following curve 
\begin{eqnarray}
\label{cornu}
  \alpha(s) 
		&=& \Bigg( 
						\frac{\pi}{\mu^2} 
						\left(\mu s \Ai(\mu s) \Bi(\mu s) 
						      - \Ai^{\prime}(\mu s) \Bi^{\prime}(\mu s)\right), \nonumber\\
		& & \quad 
		    \frac{1}{2\mu^3}
		    \left( \pi^2 
						\left( \mu s \Bi^2(\mu s) - \Bi^{\prime 2}(\mu s)\right) 
		           - \mu^3 s \Ai^2(\mu s) + \mu^2 \Ai^{\prime 2}(\mu s)\right), \nonumber\\
		& & \quad 
		    \frac{1}{2\mu^3}
		    \left( \pi^2 
						\left( \mu s \Bi^2(\mu s) - \Bi^{\prime 2}(\mu s)\right) 
		           + \mu^3 s \Ai^2(\mu s) - \mu^2 \Ai^{\prime 2}(\mu s)\right) \Bigg), 
\end{eqnarray}
if the the initial condition at $s_0=0$ is 
\begin{eqnarray*}
  \alpha(0) = \frac{1}{2 \sqrt[3]{9}\, \mu^3\, \Gamma^2(\frac13)}
		      \left(2 \sqrt{3}\, \mu\, \pi, \mu^2 - 3\pi^2, \null- \mu^2 - 3\pi^2
														\right).
\end{eqnarray*}

Thus, we can formulate the following:

\begin{corollary}
Null curves in $\mathbb E^3_1$ for which $\tau=-2\lambda s$, $\lambda=\const.\neq0$, are given by the formula (\ref{cornu}) with $\mu=\sqrt[3]{\lambda}$ up to the pseudo-arc parameter changies, up to the orientation of the curve, and up to the isometries of the space. 
\end{corollary}


\end{document}